\documentclass[11pt]{amsart}
\usepackage{mathrsfs,amsthm,amssymb,amsmath,url}

\theoremstyle{plain}
    \newtheorem{thm}{Theorem}
    
    \newtheorem{lem}[thm]{Lemma}
    \newtheorem{prop}[thm]{Proposition}
    
    \newtheorem{fact}[thm]{Fact}
    
    \newtheorem{prob}[thm]{Problem}

\theoremstyle{definition}
    \newtheorem{defn}[thm]{Definition}
    
\theoremstyle{remark}
    \newtheorem{rem}[thm]{Remark}

\newcommand{\nin}{\notin}

\newcommand{\od}{\vee}

\newcommand{\sm}{{\setminus}}

\newcommand{\cl}[1]{\langle #1 \rangle}
 
\DeclareMathOperator{\Cl}{Cl} \DeclareMathOperator{\id}{id}

\renewcommand{\O}{{\mathscr O}}
\newcommand{\On}{{\mathscr O}^{(n)}}

\newcommand{\Oo}{{\mathscr O}^{(1)}}

\newcommand{\C}{{\mathscr C}}

\newcommand{\F}{{\mathscr F}}
\newcommand{\J}{{\mathscr J}}
\newcommand{\I}{{\mathscr I}}
\newcommand{\A}{{\mathscr A}}

\newcommand{\G}{{\mathscr G}}
\newcommand{\M}{{\mathscr M}}

\renewcommand{\S}{{\mathscr S}}

\renewcommand{\L}{{\frak L}}
\newcommand{\Pos}{{\frak P}}

\newcommand{\un}{^{(n)}}

\newcommand{\uo}{^{(1)}}
\newcommand{\ut}{^{(2)}}

\newcommand{\ow}{\text{otherwise}}

\newcommand{\mpq}{m_p^{q_1,q_2}}
\newcommand{\phiqo}{\phi_{q_1}}
\newcommand{\phiqt}{\phi_{q_2}}

\author[M.\,Pinsker]{Michael Pinsker}
\address{Algebra\\TU Wien\\Wiedner Hauptstra\ss e 8-10/104\\A-1040 Wien, Austria}
\email{marula@gmx.at}\urladdr{http://dmg.tuwien.ac.at/pinsker/}
\thanks{The author is grateful for support through project P17812 of the Austrian Science Fund}
\title[Sublattices of the clone lattice]
    {
        Algebraic lattices are complete sublattices of the clone lattice over
        an infinite set
    }



\begin{document}

    \begin{abstract}
        The clone lattice $\Cl(X)$ over an infinite set $X$ is a
        complete algebraic lattice with $2^{|X|}$ compact elements.
        We show that every algebraic lattice with
        at most $2^{|X|}$ compact elements
        is a complete sublattice of $\Cl(X)$.
    \end{abstract}

   \maketitle


\section{How complicated is the clone lattice?}
    Fix a base set $X$ and denote for all $n\geq 1$ the set of all
    $n$-ary operations on $X$ by $\On$. Then $\O=\bigcup_{n\geq
    1}\On$ is the set of all functions on $X$ which have finite arity. A set of
    finitary functions $\C\subseteq\O$ is called a \emph{clone}
    iff it is closed under composition and contains all projections,
    i.e. for all $1\leq i\leq n$ the function $\pi^n_i$ satisfying
    $\pi^n_i(x_1,\ldots,x_n)=x_i$. The set of all clones over $X$
    forms a complete algebraic lattice $\Cl(X)$ with respect to
    inclusion. This lattice is countably infinite and completely
    known if $|X|=2$ by a result of Post's \cite{Pos41}; however,
    describing the clone lattice completely for larger $X$ is
    believed impossible.

    Several known results suggest this: To begin with, $\Cl(X)$ is large; it is of size continuum if $X$ is finite
    and has at least three elements, and $|\Cl(X)|=2^{2^{|X|}}$ if $X$ is
    infinite. Then, the clone lattice does not satisfy any
    non-trivial lattice identity if $|X|\geq  3$ \cite{Bul93}; it does
    not satisfy any quasi-identity if $|X|\geq 4$ \cite{Bul94}. Also, if
    $|X|\geq 4$, then every countable product of finite lattices
    is a sublattice of $\Cl(X)$ \cite{Bul94}. As for an example on infinite
    $X$, every completely distributive lattice having not more
    than $2^{|X|}$ compact elements is a subinterval of a monoidal interval of $\Cl(X)$
    \cite{Pin043} (a \emph{monoidal interval} being an interval of clones which have the same unary functions).

    We are interested in which lattices can be
    embedded into the clone lattice over an infinite set. Assume
    henceforth $X$ to be infinite. The compact elements of $\Cl(X)$ are easily seen to be exactly the
    clones which are generated by a finite number of functions.
    Since $|\O|=2^{|X|}$, this implies that $\Cl(X)$ has at most
    $2^{|X|}$ compact elements, and it is readily verified that there really exist $2^{|X|}$ compact
    elements. We are going to prove that $\Cl(X)$ is in some sense
    the most complicated algebraic lattice with this property.

    \begin{thm}
        Let $X$ be infinite. Then every algebraic lattice with at most
        $2^{|X|}$ compact
        elements can be completely embedded into
        $\Cl(X)$.
    \end{thm}

    We remark that the corresponding statement does not hold on finite
    $X$: There, $\Cl(X)$ has countably infinitely many compact (finitely generated)
    elements, but as has been proven in \cite{Bul01}, the
    countably infinite lattice $M_\omega$ does not embed into the
    clone lattice over any finite set.

    \subsection{Notation}
    We denote the unary projection $\pi^1_1$ by the somewhat
    simpler symbol $\id$, and use $\J$ for the set of projections on $X$. If $\F\subseteq\Oo$, then we write
    $\cl{\F}$ for the clone generated by $\F$. Three lattices will appear
    in the proof, the clone lattice $\Cl(X)$, the lattice $\L$ to be embedded
    into the clone lattice, and the lattice of join-semilattice ideals of compact
    elements of $\L$: For all of them, we
    use the symbols $\wedge,\vee,\bigwedge,\bigvee$ with their
    standard meanings, and confusion shall be carefully avoided. If
    $\Phi\subseteq\Oo$ is a set of unary operations, then $\Phi^*$
    will stand for all those functions which arise from functions
    of $\Phi$ by the addition of any finite number of dummy variables. Such functions
    will remain \emph{essentially unary}, i.e. although possibly non-unary they depend on only
    one variable, as opposed to \emph{essentially at least binary}
    functions, which are functions that depend on at least two of
    their variables.
\section{Proof of the main theorem}

    Let $\L$ be the lattice to be embedded into $\Cl(X)$ and
    denote by $\Pos$ the set of all compact elements of $\L$. Then $\Pos$ is a
    join-semilattice (cf. the textbook \cite{Gra78}). By an \emph{ideal} $I\subseteq\Pos$ we mean a lower subset of $\Pos$ closed under joins. The set of all
    ideals of $\Pos$ is a complete algebraic lattice, and in fact

    \begin{fact}
        $\L$ is isomorphic to the lattice of ideals of $\Pos$.
    \end{fact}

    We are going to assign a clone $\C_I$ to every ideal $I\subseteq \Pos$
    in such a way that the resulting mapping is a complete
    embedding of $\L$ into $\Cl(X)$. Fix four elements $0,1,2,4\in
    X$ and set $A=X\sm\{0,1,2,4\}$. Let $\A=(A_p)_{p\in \Pos}$ be a family of subsets of $A$ indexed by
    the elements of $\Pos$ having the following property: Whenever
    $A_p,A_{q_1},\ldots,A_{q_k}\in\A$ and $p\neq q_i$ for all
    $1\leq i\leq k$, then $A_p\nsubseteq A_{q_1}\cup\ldots\cup
    A_{q_k}$. Such a family exists: For example, there exist
    \emph{independent} families of size $2^{|X|}$, where a family
    $\F$ of subsets of $A$ is called independent iff for all finite
    disjoint $\F_1,\F_2\subseteq\F$
    $$
        \bigcap\{F: F\in\F_1\}\cap\bigcap \{A\setminus F:F\in \F_2\}\neq\emptyset.
    $$
    See the textbook \cite[Lemma 7.7]{Jec02}. If $|X|=\aleph_0$, then one could also take $\A$ to be \emph{almost
    disjoint}, meaning that all members of $\A$ are infinite and
    the intersection of any two distinct sets from $\A$ is
    finite (cf. \cite[Lemma 9.21]{Jec02}).\\

    Define for all $p\in\Pos$ a unary function $\phi_p$ by
    $$
        \phi_p(x)=
        \begin{cases}
            0&,x\in A\sm A_p\\
            1&,x\in A_p\\
            2&,x=2\\
            4&,x\in\{0,1,4\};
        \end{cases}
    $$
    so on $A$, $\phi_p$ is the characteristic function of $A_p$. Set
    $\Phi=\{\phi_p:p\in\Pos\}$.
    Now define for all $p,q_1,q_2\in\Pos$
    with $p\leq q_1\vee q_2$ a ternary function $m_p^{q_1,q_1}$ by
    $$
        m_p^{q_1,q_1}(x,y,z)=
        \begin{cases}
            \phi_p(x)&  ,y=\phi_{q_1}(x)\wedge z=\phi_{q_2}(x)\\
            2&          ,(x=2\vee y=2\vee z=2)\wedge(y\nin\{1,4\})\wedge (z\nin\{1,4\})\\
            4&          ,\ow.
            \end{cases}
    $$

    The function is well-defined: We only have to check that there is no conflict
    between the conditions for $\mpq(x,y,z)$ to
    yield $\phi_p(x)$ and $2$, respectively. If both conditions
    are satisfied, then one of the components of the tuple  the
    tuple $(x,y,z)$ equals $2$; since
    $y=\phi_{q_1}(x)$ and $z=\phiqt(x)$, this implies $x=y=z=2$,
    making the function value $\mpq(x,y,z)=2=\phi_p(x)$ unique.

    We write $\M=\{m_p^{q_1,q_1}:
    p,q_1,q_2\in\Pos\wedge p\leq q_1\vee q_2\}$ and
    $\C=\cl{\Phi\cup\M}$. The following lemma follows easily by induction over terms in $\C$.

    \begin{lem}\label{LEM:rangeInA}
        The only functions in $\C$ which take values in $A$ are the projections.
    \end{lem}

    \begin{defn}
        We call a function $f\in\Oo$ \emph{distracted} iff there
        exists $a\in A$ such that $f(a)\in\{2,4\}$.
    \end{defn}

    \begin{lem}\label{LEM:insertingDistracted}
        Let $t\in\C\un$ and $t_1,\ldots,t_n\in\Oo$. If $t$ depends
        on its $i$-th variable, where $1\leq i\leq n$, and if $t_i$ is distracted, then
        $t(t_1,\ldots,t_n)$ is distracted.
    \end{lem}
    \begin{proof}
        We use induction over terms in $\C$. To start with, let
        $t\in\J\cup\Phi\cup\M$. There is nothing to show if $t$ is a projection. If $t\in\Phi$ and $t_1\in\Oo$ is
        distracted, then there exists $a\in A$ such that
        $t_1(a)\in\{2,4\}$, so $t(t_1(a))\in\{2,4\}$
        and $t(t_1)$ is distracted. If $t=\mpq\in\M$ and $t_i$ is
        distracted for some $i\in\{1,2,3\}$, then
        $t_i(a)\in\{2,4\}$ for some $a\in A$ implies that
        $\mpq(t_1,t_2,t_3)(a)\in\{2,4\}$: Indeed, if
        $\mpq(t_1,t_2,t_3)(a)\in\{0,1\}$, then the definition of $\mpq$ would allow us to conclude $t_1(a)\in A$
        and $t_2(a)=\phiqo(t_1(a))\in\{0,1\}$ and $t_3(a)=\phiqt(t_1(a))\in\{0,1\}$,
        which is clearly impossible as
        $t_i(a)\in\{2,4\}$.\\
        For the induction step, assume that $t=f(s_1,\ldots,s_m)$,
        where $f\in\J\cup\Phi\cup\M$ and $s_j$ satisfies the induction
        hypothesis, $1\leq j\leq m$. Now there exists $1\leq j\leq m$ such that $f$ depends on its $j$-th variable and $s_j$
        depends on its $i$-th variable. By induction hypothesis $s_j(t_1,\ldots,t_n)$
        is distracted and so is
        $f(s_1(t_1,\ldots,t_n),\ldots,s_m(t_1,\ldots,t_n))$, by
        the same proof as for the induction beginning.
    \end{proof}

    \begin{lem}\label{LEM:composition}
        Let $m_p^{q_1,q_2}\in \M$ and $t_1,t_2,t_3\in\Phi\cup\{\id\}$. Then $f=m_p^{q_1,q_2}(t_1,t_2,t_3)$ is
        distracted unless $t_1=\id$, $t_2=\phi_{q_1}$, and
        $t_3=\phi_{q_2}$. In the latter case we have $f=\phi_p$.
    \end{lem}
    \begin{proof}
        If $t_2=\id$ or $t_3=\id$, then $f(a)\in\{2,4\}$ for all
        $a\in A$, since $m_p^{q_1,q_2}$ can yield $0$ or $1$ only if its second and third argument
        is in the range of a function in $\Phi$;
        hence $f$ is distracted in that case. Assume henceforth
        $t_2,t_3\in\Phi$ and write $t_2=\phi_r$ and $t_3=\phi_s$, where
        $r,s\in\Pos$.\\
        If $t_1=\id$, then $f$ yields $4$ on the symmetric differences $A_{q_1}\Delta
        A_r$ and $A_{q_2}\Delta A_s$ by the very definition of
        $m_p^{q_1,q_2}$. Hence $f$ is distracted unless those sets
        are empty, i.e. $s=q_1$ and $r=q_2$; in the latter case we
        have $f=\phi_p$ as asserted.\\
        If $t_1=\phi_l\in\Phi$, then $\mpq(\phi_l,\phi_{r},\phi_{s})$ yields by definition
        either $2,4$, or an element of the form
        $\phi_p(\phi_l(x))\in\{2,4\}$, so $f$ is distracted.
    \end{proof}

    \begin{lem}\label{LEM:allUnaryDistracted}
        All $t\in\C\uo\sm(\Phi\cup \{\id\})$ are distracted.
    \end{lem}
    \begin{proof}
        We prove this by induction over terms in $\C$. The
        beginning is trivial since there are no unary functions in the
        generating set $\J\cup\Phi\cup\M$ of $\C$ except those from
        $\Phi\cup\{\id\}$.\\
        For the induction step, assume that $t=f(t_1,\ldots,t_n)$,
        where $f\in\J\cup\Phi\cup\M$ and $t_i$ satisfies the induction
        hypothesis, for all $1\leq i\leq n$. The case $f\in\J$ is trivial. If $f\in\Phi$ and $t_1\neq\id$,
        then $t_1$ takes only values outside $A$ by Lemma
        \ref{LEM:rangeInA}, so $f(t_1)$ takes only values in $\{2,4\}$ and is distracted.
        The other possibility is that $f\in\M$, so write
        $t=m_p^{q_1,q_2}(t_1,t_2,t_3)$.
        If any of the $t_i$ is distracted then so is $t$, by Lemma \ref{LEM:insertingDistracted}.
        We may therefore assume that the
        $t_i$ are not distracted and hence elements of
        $\Phi\cup\{\id\}$. But then Lemma \ref{LEM:composition}
        tells us that $t$, not being an element of
        $\Phi\cup\{\id\}$ by assumption, must be distracted.
    \end{proof}

    \begin{defn}
        We say that $t\in\C\un$ is \emph{unspoilt} iff there exist
        $t_1,\ldots, t_n\in\C\uo$ such that
        $t(t_1,\ldots,t_n)\in\Phi$. Otherwise we call $t$ \emph{spoilt}.
    \end{defn}

    \begin{rem}
        By Lemmas \ref{LEM:insertingDistracted} and
        \ref{LEM:allUnaryDistracted}, $t_i$ must be in
        $\Phi\cup\{\id\}$ if $t$ depends on its $i$-th variable,
        for all $1\leq i\leq n$.
    \end{rem}
    \begin{rem}
        An easy induction using Lemmas
        \ref{LEM:insertingDistracted} and
        \ref{LEM:composition}
        shows that $t_i$ is uniquely determined if $t$ depends on its $i$-th variable, for all $1\leq i\leq n$.
    \end{rem}
    \begin{rem}
        By Lemmas
        \ref{LEM:insertingDistracted} and \ref{LEM:allUnaryDistracted}, a unary $t\in\C\uo$ is distracted
        iff it is spoilt.
    \end{rem}

    \begin{lem}\label{LEM:inserting2}
        Let $t\in\C\un$ be unspoilt, and assume it depends on its first variable.
        Then $t(2,x_2,\ldots,x_n)\in\{2,4\}$
        for all $x_2,\ldots,x_n\in X$.
    \end{lem}
    \begin{proof}
        We use induction over the complexity of $t$. The lemma is trivial if
        $t\in\J\cup\Phi\cup\M$. For the induction step, since the range of $\phi_p(t_1)$ is contained in
        $\{2,4\}$ and since therefore $\phi_p(t_1)$ is spoilt for all $\phi_p\in\Phi$ and all $t_1\in\C\sm\J$,
        we may assume
        $t=\mpq(t_1,t_2,t_3)$, where $t_i$ satisfies the induction
        hypothesis, $1\leq i\leq 3$. Now one of the $t_i$ must
        depend on its first variable, implying
        $t_i(2,x_2,\ldots,x_n)\in\{2,4\}$ by induction hypothesis.
        Hence, $\mpq(t_1,t_2,t_3)(2,x_2,\ldots,x_n)\in\{2,4\}$ by
        definition of $\mpq$.
    \end{proof}

    Let $t(x,y)\in\C\ut$, and consider a concrete representation $r=r(t)$ of $t$ as a term over the generating set
    $\J\cup\Phi\cup\M$ of $\C$. In the following, we write such
    representations without the use of projections, using the
    variables $x,y$ instead: For example, we write $\mpq(x,y,y)$ instead
    of $\mpq(\pi^2_1,\pi^2_2,\pi^2_2)$. This is no loss of generality and only avoids
    unnecessary usage of the projections, as in
    $\pi^2_1(\pi^2_2,\phi_p(\pi^2_1))$ (equivalently, we could demand the projections to appear only as innermost arguments in the
    representation). We say that a subterm $s$ of $r$ is a
    \emph{leaf} of $r$ iff it involves exactly one function
    symbol from $\Phi\cup\M$. For example, the leaves of
    $$
        \mpq(m_u^{v_1,v_2}(x,\phi_l(y),\phi_r(x)),\phi_d(y),m_g^{h_1,h_2}(x,x,x))
    $$
    are $\phi_l(y),\phi_r(x),\phi_d(y),$ and
    $m_g^{h_1,h_2}(x,x,x)$. Thinking of $r$ as a tree in which the variables are not represented by an own node,
    the leaves of $r$ are really exactly the leaves of the
    tree.\\
    We call the representation $r(t)$ \emph{reduced} iff it
    has no subterms of the form $\mpq(x,\phiqo(x),\phiqt(x))$.
    Such subterms can be replaced by $\phi_p(x)$ by virtue of Lemma \ref{LEM:composition}, so every term $t$ has
    a reduced representation. We are only interested in
    representations of unspoilt functions that depend on both variables, so all unary subterms
    of any representation
    correspond to elements of $\Phi$, by Lemmas
    \ref{LEM:insertingDistracted} and
    \ref{LEM:allUnaryDistracted}; working with reduced terms means that we demand those unary subterms
    to be represented by only one function symbol.\\
    Let $r(t)$ be reduced. We set $Leaf(r)$ to consist of all
    leaves of $r(t)$. Note that $Leaf(r)$ depends on the
    representation of the function $t$.

    \begin{lem}\label{LEM:when4}
        Let $r(x,y)$ be a reduced representation of a binary function in $\C$ that is unspoilt and depends on both of its
        variables. Let $a\in A$. Then $r(2,a)=4$ iff
        $a\in\bigcup\{A_v: \phi_v(y)\in Leaf(r)\}$.
    \end{lem}
    \begin{proof}
        We use induction over the complexity of $r$. The beginning is trivial as there are no binary functions
        depending on both variables in the generating set of $\C$.
        For the induction step, write $r=f(r_1,\ldots,r_n)$, where
        $f\in\Phi\cup\M$, and where $r_i$ satisfies the induction
        hypothesis, $1\leq i\leq n$. If $f\in\Phi$, then using Lemma \ref{LEM:rangeInA}
        it is readily verified that $f(r_1)$ is
        spoilt unless $r_1$ is a projection, in which case
        $r\in\Phi^*$, contradicting that $r$ depends on both
        variables. Assume henceforth that $f=\mpq\in\M$.\\
        Observe that all $r_i$ must be unspoilt, for otherwise $r$
        would be spoilt as well by Lemmas \ref{LEM:insertingDistracted} and
        \ref{LEM:allUnaryDistracted}. Since $r$ is unspoilt, there
        exist $s_1,s_2\in\C\uo$ such that
        $\mpq(r_1(s_1,s_2),\ldots,r_3(s_1,s_2))\in\Phi$. By Lemmas
        \ref{LEM:insertingDistracted},
        \ref{LEM:composition} and \ref{LEM:allUnaryDistracted}, this is only possible if
        $r_1(s_1,s_2)$ is the identity, which together with Lemma \ref{LEM:rangeInA} implies that $r_1$
        is a projection. Suppose that $r_2=r_1=\pi^2_i$, where $i\in\{1,2\}$.
        Then $r(s_1,s_2)=\mpq(s_i,s_i,r_3(s_1,s_2))\in\Phi$ and
        Lemma \ref{LEM:composition} implies that the first argument in $\mpq$ must be the identity, while the second
        must equal $\phi_{q_1}$, an obvious contradiction. The same contradiction occurs assuming $r_3=r_1$, and hence
        we have $r_i\neq r_1$, $i=2,3$. We now distinguish six
        cases.\\
        Assume first that $r_2,r_3\in\J$. Then $r=\mpq(x,y,y)$ or
        $r=\mpq(y,x,x)$. In either case we have $r(2,a)=2\neq 4$, in accordance with our assertion as $r$ does not
        have any leaves of the form $\phi_v(y)$.\\
        Consider the case where $r_2\in\J$ and $r_3\in\Phi^*$ (by symmetry,
        this also treats the case $r_3\in\J$ and $r_2\in\Phi^*$).
        Keeping Lemma \ref{LEM:composition} and the fact that $r$ depends on both variables in mind we conclude
        that $r=\mpq(x,y,\phiqt(x))$ or $r=\mpq(y,x,\phiqt(y))$ or
        $r=\mpq(x,y,\phiqt(y))$ or $r=\mpq(y,x,\phiqt(x))$. The
        latter two possibilities, however, are spoilt as
        substitution of $\phiqo$ for $y$ and $x$, respectively, yields a
        distracted third argument $\phiqt(\phiqo)$ of $\mpq$. The first possibility
        gives us $r(2,a)=\mpq(2,a,2)=2\neq 4$, in
        accordance with our assertion. Finally, for the second
        term we have $r(2,a)=\mpq(a,2,\phiqt(a))$, which equals
        $4$ iff $\phiqt(a)\nin\{1,4\}$ iff $a\in A_{q_2}$.\\
        Now assume that $r_2\in\J$ and $r_3\nin\J\cup\Phi^*$. Then $r_3$ depends on both of its variables by Lemma
        \ref{LEM:allUnaryDistracted}, and therefore satisfies the assertion of this lemma by induction hypothesis.
        By Lemma
        \ref{LEM:inserting2} we have that $r(2,a)=4$ iff $r(2,a)\neq
        2$; the definition of $\mpq$ tells us that this is the
        case iff $2\nin\{r_1(2,a),r_2(2,a),r_3(2,a)\}$ or
        $r_2(2,a)\in\{1,4\}$ or $r_3(2,a)\in\{1,4\}$. Now $r_3(2,a)\in\{2,4\}$ by Lemma
        \ref{LEM:inserting2}, and $r_2(2,a)\in\{2,a\}$ since $r_2$ is a
        projection. Thus, $r(2,a)=4$ iff
        $r_3(2,a)=4$, which by induction hypothesis is the case
        iff $a\in \bigcup\{A_v: \phi_v(y)\in Leaf(r_3)\}$. Since
        the leaves of $r_3$ are the just the leaves of $r$ we are
        done.\\
        Next say that $r_2\in\Phi^*$ and $r_3\nin\J\cup\Phi^*$.
        We have $r(2,a)=4$ iff $r(2,a)\neq 2$, which happens iff $2\nin\{r_1(2,a),r_2(2,a),r_3(2,a)\}$ or
        $r_2(2,a)\in\{1,4\}$ or $r_3(2,a)\in\{1,4\}$. Again,
        $r_3(2,a)\in\{2,4\}$ by Lemma \ref{LEM:inserting2}, and $r_2(2,a)\in\{0,1,2\}$ as
        $r_2\in\Phi^*$,
        implying $r(2,a)=4$ iff $r_2(2,a)=1$ or $r_3(2,a)=4$. Now
        if $r_2(x,y)=\phi_{q_1}(x)$, then $r_2(2,a)=2$ and so $r(2,a)=4$ iff
        $r_3(2,a)=4$ iff $a\in \bigcup\{A_v: \phi_v(y)\in
        Leaf(r_3)\}$ by induction hypothesis. This is in
        accordance with our assertion since then $\phi_v(y)\in Leaf(r_3)$ iff $\phi_v(y)\in
        Leaf(r)$. If on the other hand $r_2(x,y)=\phi_{q_1}(y)$,
        then $r_2(2,a)=1$ iff $a\in A_{q_1}$, and hence $r(2,a)=4$ iff
        $a\in A_{q_1}\cup \bigcup\{A_v: \phi_v(y)\in
        Leaf(r_3)\}$; this is the case iff $a\in \bigcup\{A_v: \phi_v(y)\in
        Leaf(r)\}$.\\
        If $r_2,r_3\in\Phi^*$, then up to symmetry
        $r=\mpq(x,\phi_{q_1}(x),\phi_{q_2}(y))$ or
        $r=\mpq(x,\phi_{q_1}(y),\phi_{q_2}(y))$ or
        $r=\mpq(y,\phi_{q_1}(x),\phi_{q_2}(x))$ or
        $r=\mpq(y,\phi_{q_1}(x),\phi_{q_2}(y))$. Therefore
        $r(2,a)=4$ iff $a\in A_{q_2}$ in the first case, iff $a\in A_{q_1}\cup
        A_{q_2}$ in the second case, and iff $a\in A_{q_2}$ in the
        fourth case; in the third case, $r(2,a)=2\neq 4$.\\
        Finally, consider $r_2,r_3\nin\J\cup\Phi^*$. By Lemma
        \ref{LEM:inserting2}, $\{r_2(2,a),r_3(2,a)\}\subseteq \{2,4\}$;
        thus, $r(2,a)=4$ iff $r(2,a)\neq 2$ iff $r_2(2,a)=4$ or
        $r_3(2,a)=4$. Using the induction hypothesis, we get that
        $r(2,a)$ yields $4$ iff $a\in\bigcup\{A_v: \phi_v(y)\in Leaf(r_2)\}$ or
        $a\in\bigcup\{A_v: \phi_v(y)\in Leaf(r_3)\}$; hence,
        $r(2,a)=4$ iff $a\in\bigcup\{A_v: \phi_v(y)\in Leaf(r)\}$.
    \end{proof}

    Set $\S=\{t\in\C: t\, \text{ spoilt}\}$. Define for all $I\subseteq\Pos$ sets of functions
    $\Phi_I=\{\phi_p\in\Phi:p\in I\}$ and $\G_I=\Phi_I\cup\M\cup\S$, and a clone
    $\C_I=\cl{\G_I}$. Write $\cl{I}$ for the ideal of $\Pos$ generated by
    $I$.
    \begin{lem}\label{LEM:generatedIdeal}
        Let $p\in\Pos$ and $I\subseteq\Pos$. Then $\phi_p\in \C_I$ iff
        $p\in\cl{I}$.
    \end{lem}
    \begin{proof}
        Let $t\in\C_I$; using induction over the complexity of $t$ as a term over the generating set $\G_I$,
        we show that
        $t=\phi_p$ implies $p\in\cl{I}$. The beginning is trivial,
        since if $t\in\G_I$, then $t\in\Phi_I$ and so $p\in I$.
        For the induction step, write $t=f(t_1,\ldots,t_n)$, with
        $f\in\G_I$ and $t_i\in\C_I$ satisfying the induction hypothesis, $1\leq i\leq n$. Clearly, $f\in\S$ is impossible. $f\in\Phi_I$
        implies that $t_1$ is the identity and so $f=\phi_p$;
        hence $p\in I$. Assume therefore that $f=m_{u}^{q_1,q_2}\in\M$.
        Then $u=p$,
        $t_1=\id, t_2=\phi_{q_1}$ and $t_3=\phi_{q_2}$ by Lemmas
        \ref{LEM:insertingDistracted},
        \ref{LEM:composition} and \ref{LEM:allUnaryDistracted}. By induction hypothesis,
        $q_1,q_2\in\cl{I}$. Hence, $p\leq q_1\od q_2\in\cl{I}$.\\
        For the other direction, it is enough to show that if
        $\phi_{q_1},\phi_{q_2}\in\C_I$, then $\phi_u\in\C_I$ for
        all $u\leq q_1\od q_2$. But this is clear since
        $\phi_u=m_u^{q_1,q_2}(\id,\phi_{q_1},\phi_{q_2})\in\C_I$.
    \end{proof}

    \begin{lem}\label{LEM:supremum}
        Let $\I$ be a family of ideals of $\Pos$.
        Then $\bigvee \{\C_I:I\in\I\} =\C_{\bigvee\I}$.
    \end{lem}
    \begin{proof}
        Trivially, $\C_{\bigvee\I}$ contains all $\C_I$, where $I\in\I$,
        hence it contains $\bigvee \{\C_I:I\in\I\}$. For the other inclusion we have to show that
        $\C_{\bigvee\I}$ is contained in $\bigvee \{\C_I:I\in\I\}$; clearly, it is enough to show that
        $\Phi_{\bigvee\I}\subseteq\bigvee \{\C_I:I\in\I\}$.
        Indeed, if $\phi_p\in\Phi_{\bigvee\I}$, then $p\in
        \bigvee\I$. Since $\bigvee\I=\cl{\bigcup\I}$, the
        preceding lemma implies $\phi_p\in\C_{\bigcup\I}$. Now it
        is enough to observe that $\C_{\bigcup\I}$ equals
        $\cl{\bigcup \{\C_I:I\in\I\}}$, which is exactly $\bigvee
        \{\C_I:I\in\I\}$.
    \end{proof}

    \begin{lem}\label{LEM:infimum}
        Let $\I$ be a family of ideals of $\Pos$. Then
        $\bigwedge \{\C_I:I\in\I\}=\C_{\bigwedge\I}$.
    \end{lem}
    \begin{proof}
        $\C_{\bigwedge\I}$ is a subclone of all $\C_I$, where $I\in\I$, so
        trivially $\C_{\bigwedge\I}\subseteq \bigwedge \{\C_I:I\in\I\}$. For the
        other direction, let $t\in\bigwedge \{\C_I:I\in\I\}=\bigcap \{\C_I:I\in\I\}$. If $t$ is spoilt,
        then $t\in \C_{\bigwedge\I}$ by definition, so assume that $t$
        is unspoilt. If $t$ is essentially unary, then $t$ is a
        projection or an element of $\Phi^*$, by Lemma
        \ref{LEM:allUnaryDistracted}. In the latter case,
        $t\in\bigcap\{\Phi_I^*:I\in\I\}$ by Lemma \ref{LEM:generatedIdeal}, so $t\in \C_{\bigcap\I}=\C_{\bigwedge\I}$. So let
        $t$ be essentially at least binary, and assume without loss of generality that it depends on all of its variables.
        Because $t$ is unspoilt, there exist
        $t_1,\ldots,t_n\in\Phi\cup\{id\}$ such that
        $t(t_1,\ldots,t_n)\in\Phi$. Set
        $s_i(x,y)=t(t_1(x),\ldots,t_{i-1}(x),y,t_{i+1}(x),\ldots,t_n(x))$,
        for all $1\leq i\leq n$. Obviously, all $s_i$ are
        unspoilt. They also depend on both variables: Indeed, let
        without loss of generality $i=1$. Then
        $s_1(2,t_1(a))=t(t_1(a),2,\ldots,2)\in\{2,4\}$ by Lemma \ref{LEM:inserting2} but
        $s_1(a,t_1(a))=t(t_1,\ldots,t_n)(a)\in\{0,1\}$ for all $a\in A$,
        so $s_1$ depends on the first variable. For the second variable, observe that
        $s_1(a,2)=t(2,t_2(a),\ldots,t_n(a))\in\{2,4\}$, so
        $s_1(a,t_1(a))\neq s_1(a,2)$.\\
        Assume that $t$ is represented as a reduced term. The $s_i$ might not be
        reduced: For example, $t$ could have a subterm like
        $\mpq(x_2,\phi_{q_1}(x_3),x_4)$, which becomes
        $\mpq(x,\phi_{q_1}(x),\phi_{q_2}(x))$ when we substitute
        $x_2=x_3=x$ and $x_4=\phi_{q_2}(x)$ upon building, say, $s_1$. However, such
        redundancies will occur only for the variable
        $x$. Thus, when simplifying $s_i$ to a reduced term according to the equation $\mpq(x,\phiqo(x),\phiqt(x))=\phi_p(x)$, the leaves of the form $\phi_p(y)$, which were
        originally (that is, in $t$) leaves of the form
        $\phi_p(x_i)$, do not change. Therefore, $\phi_p(y)$ is a
        leaf of the new reduced $s_i$ iff $\phi_p(x_i)$ is a leaf of $t$.\\
        By Lemma \ref{LEM:when4}, for all $1\leq i\leq n$ and for all $a\in A$ we have that $s_i(2,a)=4$ iff
        $a\in\bigcup\{A_v: \phi_v(y)\in
        Leaf(s_i)\}$. This is the case iff $a\in\bigcup\{A_v: \phi_v(x_i)\in
        Leaf(t)\}$. Therefore, there exists $1\leq i\leq n$ with
        $s_i(2,a)=4$ iff $a\in\bigcup\{A_v: \exists i\, (\phi_v(x_i)\in
        Leaf(t))\}$. Pick arbitrary $I,J\in\I$ and
        consider two reduced representations $t_I, t_J$ of $t$, where $t_I$ is a term over $\G_I$ and
        $t_J$ one over $\G_J$. Then, since whether or not
        $s_i(2,a)=4$ does not depend on the representation,
        $$
            \bigcup\{A_v: \exists i\, (\phi_v(x_i)\in Leaf(t_I))\}=\bigcup\{A_v: \exists i\, (\phi_v(x_i)\in
            Leaf(t_J))\}.
        $$
        Because $A_v\nsubseteq A_{q_1}\cup\ldots \cup A_{q_k}$ whenever
        $q_i\neq v$, $1\leq i\leq k$, we conclude
        $$
            \{v: \exists i\, (\phi_v(x_i)\in Leaf(t_I))\}=\{v: \exists i\, (\phi_v(x_i)\in
            Leaf(t_J))\}.
        $$
        Thus, the latter set is a subset of both $I$ and $J$, implying
        that $t_I$ actually involves only functions from $\G_{I\cap J}$ as leaves. Since $J$ was arbitrary, we may conclude that
        the term $t_I$ uses only functions from $\G_{\bigwedge \I}$ as leaves. Because functions from $\Phi$
        can appear only as leaves in an unspoilt term ($\phi_v(f)$ is spoilt for all $\phi_v\in\Phi$ and all $f\in\C$ unless $f$ is a projection), this means that $t_I$ contains only functions
        from $\G_{\bigwedge \I}$. Hence, $t\in\C_{\bigwedge \I}$.
    \end{proof}

    \begin{prop}\label{PROP:embedding}
        The mapping assigning $\C_I$ to every ideal $I\subseteq \Pos$
        is a complete lattice embedding of $\L$
        into $\Cl(X)$.
    \end{prop}
    \begin{proof}
         The function is injective by Lemma
         \ref{LEM:generatedIdeal} and preserves arbitrary suprema and infima by Lemmas
         \ref{LEM:supremum} and \ref{LEM:infimum}.
    \end{proof}

\section{Concluding remarks and outlook}
    The only place where we used the infinity of the base set $X$
    is when we claim the existence of a family $\A$ which is as large as $\Pos$ and has
    the property that whenever
    $A_p,A_{q_1},\ldots,A_{q_k}\in\A$ and $p\neq q_i$ for all
    $1\leq i\leq k$, then $A_p\nsubseteq A_{q_1}\cup\ldots\cup
    A_{q_k}$.
    Therefore surprisingly, the same proof works to show that every finite lattice $\L$ is a sublattice of the
    clone lattice over a finite $X$ for some $X$ large enough ($|X|\geq |\L|+4$ suffices).
    However, as mentioned in the introduction, much better results already exist for finite $X$.

    Answering the following question would be a next interesting
    step in answering the question of how complicated the clone
    lattice is.

    \begin{prob}
        Is every algebraic lattice with at most $2^{|X|}$ compact elements an interval of $\Cl(X)$?
    \end{prob}

    \bibliographystyle{alpha}

\end{document}